\documentclass{elsarticle}
\usepackage{hyperref}
\usepackage{mathrsfs}
\let\memoldbibsection\bibsection
\let\bibsection\relax
\usepackage{amsrefs}         
\let\bibsection\memoldbibsection
\usepackage{amsmath}

\def\cl{\mathop{\rm cl}\nolimits}
\def\CO{\mathop{\rm CO}\nolimits}

\def\int{\mathop{\rm int}\nolimits}

\def\supp{\mathop{\rm supp}\nolimits}

\usepackage{commath}
\usepackage{xcolor}

\usepackage{tikz, tikz-cd}

\makeatletter
\def\ps@pprintTitle{%
  \let\@oddhead\@empty
  \let\@evenhead\@empty
  \let\@oddfoot\@empty
  \let\@evenfoot\@oddfoot
}
\makeatother

\usepackage{amsfonts}
\usepackage{amssymb}
\newtheorem{thm}{Theorem}
\newtheorem{lem}[thm]{Lemma}

\newtheorem{pro}{Proposition}
\newdefinition{rmk}{Remark}
\newtheorem{cor}{Corollary}
\newtheorem{question}{Question}

\newproof{pf}{Proof}

\begin{document}

\begin{frontmatter}

\title{Remarks on SHD spaces and more divergence properties}

\author[1]{Carlos David Jiménez-Flores\fnref{fn1}}
\ead{carlosjf@ciencias.unam.mx}

\author[1]{Alejandro Ríos-Herrejón\fnref{fn1}}
\ead{chanchito@ciencias.unam.mx}

\author[2]{Alejandro Darío Rojas-Sánchez}
\ead{adrojas@up.edu.mx}

\author[3]{Artur Hideyuki Tomita\fnref{fn2}}
\ead{tomita@ime.usp.br}

\author[1]{Elmer Enrique Tovar-Acosta\fnref{fn1}}
\ead{elmer@ciencias.unam.mx}

\affiliation[1]{organization={Universidad Nacional Autónoma de México.},
            addressline={Av. Universidad 3000}, 
            city={CDMX},
            postcode={4510}, 
            country={México}}  

\affiliation[2]{organization={Universidad Panamericana},
            addressline={Augusto Rodin No. 498, Col. Insurgentes Mixcoac}, 
            city={CDMX},
            postcode={03920}, 
            country={México}}

\affiliation[3]{organization={University of São Paulo},
            addressline={Rua do Matão}, 
            city={São Paulo},
            postcode={1010}, 
            country={Brazil}}

\fntext[fn1]{The first, second and fifth authors acknowledge financial support from CONAHCyT grants no. 816193, 814282 and  829699, respectively.}

\fntext[fn2]{The fourth author acknowledges financial support from FAPESP 2021/00177-4}

\begin{abstract}
The class of SHD spaces was recently introduced in \cite{SHD}. The first part of this paper focuses on answering most of the questions presented in that article. For instance, we exhibit an example of a non-SHD Tychonoff space $X$ such that $\mathscr{F}[X]$, the Pixley-Roy hyperspace of $X$, $\beta X$, the Stone-\v{C}ech compactification of $X$, and $C_p(X)$, the ring of continuous functions over $X$ equipped with the topology of pointwise convergence, are SHD.

In the second part of this work we will present some variations of the SHD notion, namely, the WSHD property and the OHD property. Furthermore, we will pay special attention to the relationships between $X$ and $\mathscr{F}[X]$ regarding these new concepts.

\end{abstract}


\end{frontmatter}

\section{Remarks on SHD spaces}

All set-theoretic concepts whose definition is not included here must be assumed as in \cite{kunen2013}. Our basic reference book for topological notions is \cite{Engelking}, however, the specific material corresponding to cardinal functions should be understood as in \cite{hodel1984}. Throughout this paper $X$ will be a topological space with at least two points, $\tau_X$ will denote its topology and $\tau_X^+$ will be the family $\tau_X \setminus \{\emptyset\}$.

Recall that if $\{X_\alpha : \alpha<\kappa\}$ is a familiy of topological spaces, then the collection $$\mathscr{B}\left(\prod_{\alpha<\kappa} X_\alpha\right) := \left\{\bigcap_{\alpha \in F} \pi_{\alpha}^{-1}[U_\alpha] : F\in [\kappa]^{<\omega}\setminus \{\emptyset\} \wedge \forall \alpha \in F\left(U_\alpha \in \tau_{X_\alpha}^{+}\right)\right\}$$ is the {\it canonical base} for the standard product topology. If $U\in \mathscr{B}\left(\prod_{\alpha<\kappa} X_\alpha\right)$, then the {\it support of $U$} is the set $$\supp(U) := \left\{\alpha<\kappa : \pi_\alpha[U] \neq X_\alpha \right\}.$$Finally, $P(X)$ is used to denote the power set of $X$, i.e., the set of all subsets of $X$.

The class of SHD spaces appeared for the first time in \cite{SHD}. Recall that $X$ is {\it selectively highly divergent} (in short, {\it SHD}) if for every $\{U_n : n\in\omega\}\subseteq \tau_X^+$ we can pick $x_n \in U_n$ for each $n\in\omega$ in such a way that the sequence $(x_n)$ has no convergent subsequences. In this section we will answer most of the questions presented in \cite{SHD} and prove several new results. Let us begin with the following proposition.

\begin{pro}\label{generalization} Let $\kappa$ be an uncountable cardinal and $\{X_\alpha : \alpha\in\kappa\}$ be a collection of topological spaces. Suppose there is $D\subseteq \kappa$ such that $\{X_\alpha : \alpha\in D\}$ is a family of non-sequentially compact spaces. If $D\subseteq\kappa$ is uncountable then $\prod_{\alpha\in\kappa} X_\alpha$ is SHD.

\end{pro}

\begin{pf} Set $X:= \prod_{\alpha\in\kappa} X_\alpha$ and let $\{U_n : n\in\omega \}$ be a subset of $\mathscr{B}(X)$. Set $F:= \bigcup_{n\in \omega} \supp(U_n)$, fix $\beta\in D\cap (\kappa\setminus F)$ and let $(y_n)$ be a sequence in $X_\beta$ with no convergent subsequences. For each $n\in \omega$ define $x_n\in X$ by $$x_n(\alpha) := \begin{cases} y_n, & \text{if} \ \alpha = \beta, \\
\in \pi_\alpha[U_n], & \text{if} \ \alpha \in \supp(U_n), \\
\in X_\alpha, & \text{if} \ \alpha \in \kappa \setminus\left(\supp(U_n) \cup \{\beta\}\right).
  \end{cases}$$ Clearly, $x_n \in U_n$ for every $n\in \omega$. Furthermore, since $(y_n)$ has no convergent subsequences, $(x_n)$ has no convergent subsequences.

\end{pf}

\begin{cor} If $\kappa$ is an uncountable cardinal and $X$ is not sequentially compact, then $X^{\kappa}$ is SHD.

\end{cor}

Evidently, a sequentially compact space cannot be SHD. It turns out that in the case of uncountable topological powers both conditions are equivalent:

\begin{pro}\label{prop_equivalencia_potencia} If $\kappa$ is an uncountable cardinal, then $X^{\kappa}$ is SHD if and only if $X^{\kappa}$ is not sequentially compact. 

\end{pro}

\begin{pf} Suppose that $X^{\kappa}$ is not sequentially compact, let $\{U_n : n\in\omega \}$ be a subset of $\mathscr{B}\left(X^{\kappa}\right)$ and set $F:=\bigcup_{n\in\omega} \supp(U_n)$. Observe that if $I := \kappa \setminus F$, then $X^\kappa$ is homeomorphic to $X^{I}$ and thus, $X^{I}$ is not sequentially compact. Fix a sequence $(y_n)$ in $X^{I}$ with no convergent subsequences. For each $n\in \omega$ define $x_n\in X^\kappa$ by $$x_n(\alpha) := \begin{cases} y_n(\alpha), & \text{if} \ \alpha \in I, \\
\in \pi_\alpha[U_n], & \text{if} \ \alpha \in \supp(U_n), \\
\in X_\alpha, & \text{if} \ \alpha \in F \setminus \supp(U_n).
  \end{cases}$$ Certainly, $x_n \in U_n$ for every $n\in \omega$. Moreover, since $(y_n)$ has no convergent subsequences, $(x_n)$ has no convergent subsequences.

\end{pf}

Since there are no first countable SHD spaces it can be easily seen that the restriction on $\kappa$ in Proposition~\ref{prop_equivalencia_potencia} is essential. For example, if $X$ is first countable and not sequentially compact, then $X^{\omega}$ is non-SHD and non-sequentially compact.

Recall that $X$ is {\it ultraconnected} if there are no disjoint non-empty closed subsets of $X$. It is known that if $\mathfrak{s}$ stands for the {\it splitting number} (see \cite{halbeisen2012}) then the product of at least $\mathfrak{s}$ spaces which are not ultraconnected is not sequentially compact (see \cite[Lemma~4.2, p.~995]{lipparini}). Consequently, the following result holds.

\begin{cor}\label{cor_split} If $X$ is not ultraconnected and $\kappa\geq \mathfrak{s}$ is a cardinal number, then $X^{\kappa}$ is SHD.

\end{cor}

It was shown in \cite{SHD} that under the Generalized Continuum Hypothesis it happens that for every infinite cardinal $\kappa$ there is a space $X_\kappa$ that is Hausdorff, compact, zero-dimensional and SHD with $d(X_\kappa)=\kappa$ and $w(X_\kappa)=2^{\kappa}$. With the previous result available we can get the same conclusion without using additional axioms.

\begin{pro} For every infinite cardinal $\kappa$ there exists a topological space $X_\kappa$ that is Hausdorff, compact, zero-dimensional and SHD with $d(X_\kappa)=\kappa$ and $w(X_\kappa)=2^{\kappa}$.

\end{pro}

\begin{pf} Let $\alpha D(\kappa)$ stand for the one-point compactification of the discrete space of size $\kappa$. Clearly, $X_\kappa := (\alpha D(\kappa))^{2^{\kappa}}$ is Hausdorff, compact and zero-dimensional. Moreover, a standard argument with cardinal functions shows that $d(X_k)=\kappa$ and $w(X_\kappa)=2^{\kappa}$. Finally, Corollary~\ref{cor_split} implies that $X_\kappa$ is SHD since $\mathfrak{s}\leq 2^{\kappa}$.
\end{pf}

\begin{rmk} Let $X$ be a non-ultraconnected sequentially compact space. If $\mathfrak{t}$ stands for the {\it tower number} (see \cite{halbeisen2012}) and $\omega_1<\mathfrak{t}$, it turns out that $X^{\omega_1}$ is sequentially compact (see \cite[Theorem~6.9(a), p.~132]{vanDouwen1984}) and thus, $X^{\omega_1}$ is not SHD. On the other hand, $\omega_1=\mathfrak{s}$ implies that $X^{\omega_1}$ is SHD thanks to Corollary~\ref{cor_split}. Hence, the statement {\lq\lq}the $\omega_1$-th power of a non-ultraconnected sequentially compact space is SHD{\rq\rq} is independent of \textsf{ZFC}.

\end{rmk}

In the case of Cantor cubes, A. Bella and S. Spadaro proved in \cite{bella} that $D(2)^{\kappa}$ is SHD if and only if $\kappa\geq \mathfrak{s}$. Our next goal is to prove a strengthening of Bella's result.

\begin{thm}\label{thm_compacto_SHD} Let $X$ be a compact Hausdorff space and $\kappa$ be an uncountable cardinal. If $w(X) \leq \omega_1$ or $X$ is dyadic (i.e., $X$ is a continuous image of a Cantor cube) and $w(X)<\mathfrak{s}$, then the following statements are equivalent. 

\begin{enumerate}
\item $X^{\kappa}$ is SHD.
\item $X^{\kappa}$ is not sequentially compact.
\item $\kappa\geq \mathfrak{s}$.
\end{enumerate}

\end{thm}

\begin{pf} Items (1) and (2) are equivalent thanks to Proposition~\ref{prop_equivalencia_potencia}. On the other hand, Corollary~\ref{cor_split} guarantees that (3) implies (1). For our remaining implication suppose that (2) holds and let's concentrate on proving (3). If $w(X) \leq \omega_1$, then $w(X^\kappa) = \kappa$ and thus, since $X^\kappa$ is compact and not sequentially compact, it follows from \cite[Theorem~6.1, p.~129]{vanDouwen1984} that $\kappa\geq \mathfrak{s}$. 

In the case where $X$ is dyadic and $w(X)<\mathfrak{s}$, there exists a continuous surjection from $D(2)^{w(X)}$ onto $X$ (see \cite[3.12.12(b), p.~231]{Engelking}), which in turn produces a continuous surjection from $D(2)^{w(X)\cdot \kappa}$ onto $X^\kappa$. Now, since $X^\kappa$ is not sequentially compact, $D(2)^{w(X)\cdot \kappa}$ cannot be sequentially compact and therefore, Bella and Spadaro's theorem combined with Proposition~\ref{prop_equivalencia_potencia} yields the inequality $w(X)\cdot \kappa\geq \mathfrak{s}$; hence, $\kappa\geq \mathfrak{s}$.
\end{pf}

Recall that if $p$ is an element of $X$ and $\kappa$ is an infinite cardinal, then

\begin{align*} \sigma\left(X,p,\kappa\right) &:= \left\{x\in X^{\kappa} : \abs{\left\{\alpha\in\kappa : x(\alpha)\neq p\right\}}<\omega\right\} \ \text{and} \\
\Sigma\left(X,p,\kappa\right) &:= \left\{x\in X^{\kappa} : \abs{\left\{\alpha\in\kappa : x(\alpha)\neq p\right\}}\leq \omega\right\}.
\end{align*}

\begin{rmk}\label{ex_sigma_producto} Let $\kappa\geq \mathfrak{s}$ and $p\in D(2)$. It is well-known that $\Sigma\left(D(2),p,\kappa\right)$ is a sequentially compact dense subspace of $D(2)^{\kappa}$. Hence, $\Sigma\left(D(2),p,\kappa\right)$ is a non-SHD dense subspace of the SHD space $D(2)^{\kappa}$. The case $\kappa = \mathfrak{c}$ was first shown in \cite{bella}.

\end{rmk}

One of the questions discussed in \cite{SHD} was whether it was possible to establish a connection with respect to the SHD property between a Tychonoff space $X$ and $\beta X$, its Stone-\v{C}ech compactification. It turns out that $X$ does not have to be SHD even if $\beta X$ is:

\begin{thm}\label{thm_rios_tomita} If $\kappa\geq \mathfrak{s}$ is a cardinal number and $X$ is a dense subspace of $D(2)^{\kappa}$, then $\beta X$ is SHD.

\end{thm}

\begin{pf} Set $Y:= D(2)^{\kappa}$. Let $\{W_n : n\in\omega\}$ be a subset of $\tau_{\beta X}^{+}$ and pick $\{V_n : n\in\omega\} \subseteq \tau_{\beta X}^{+}$ with $V_n \subseteq \cl_{\beta X} (V_n) \subseteq W_n$ for every $n\in \omega$. Since $\{V_n \cap X : n\in \omega\}$ is a subset of $\tau_X^+$, there is a family $\{O_n : n\in \omega\}\subseteq \tau_Y^+$ such that $V_n \cap X = O_n \cap X$ whenever $n\in \omega$. For every $n\in \omega$ let $U_n\in \mathscr{B}(Y)$ satisfy $U_n \subseteq O_n$. Set $F:= \bigcup_{n\in\omega} \supp(U_n)$ and $Z:= D(2)^{\kappa\setminus F}$. Since $Z$ is not sequentially compact (see Theorem~\ref{thm_compacto_SHD}), there is a sequence $(z_n)$ in $Z$ with no convergent subsequences. For each $n\in \omega$ define $y_n\in Y$ by $$y_n(\alpha) := \begin{cases} z_n(\alpha), & \text{if} \ \alpha \in \kappa\setminus F, \\
\in \pi_\alpha[U_n], & \text{if} \ \alpha \in \supp(U_n), \\
0, & \text{if} \ \alpha \in F\setminus \supp(U_n).
  \end{cases}$$ Clearly, $y_n \in U_n$ for every $n\in \omega$.
  
Observe that since $Y$ is a compactification of $X$, there exists a continuous surjective function $\beta i : \beta X \to Y$ such that its restriction to $X$ is the identity function. We claim that $U_n$ is a subset of $\beta i[W_n]$ for every $n\in\omega$. To prove this fix $n\in \omega$ and notice three facts: first, $U_n\cap X = \beta i[U_n\cap X] \subseteq \beta i[\cl_{\beta X}(U_n\cap X)]$; second, $\cl_{\beta X}(U_n\cap X)\subseteq \cl_{\beta X}(V_n)\subseteq W_n$; and third, $\beta i \left[\cl_{\beta X}(U_n\cap X)\right]$ is a closed subset of $Y$ since $\cl_{\beta X}(U_n\cap X)$ is a compact subset of $\beta X$. Therefore, $U_n=\cl_{Y}(U_n)=\cl_{Y}(U_n\cap X)\subseteq \beta i[\cl_{\beta X}(U_n\cap X)]\subseteq \beta i[W_n]$.

For each $n\in \omega$ let $w_n \in W_n$ satisfy $\beta i (w_n) = y_n$. To finish our proof suppose that $A\in [\omega]^{\omega}$ and $w\in \beta X$ are such that $(w_n)_{n\in A}$ converges to $w$ in $\beta X$. Hence, since $\beta i : \beta X \to Y$ and the natural projection $\pi : Y \to Z$ are continuous, $(z_n)_{n\in A}$ converges to $\pi\left(\beta i (w)\right)$ in $Z$, a contradiction. Thus, $(w_n)$ has no convergent subsequences.

\end{pf}

Remark~\ref{ex_sigma_producto} and Theorem~\ref{thm_rios_tomita} yield the following corollary.

\begin{cor}\label{cor_thm_rios_tomita} If $\kappa\geq \mathfrak{s}$ and $p\in D(2)$, then $\Sigma\left(D(2),p,\kappa\right)$ is not SHD and $\beta \Sigma\left(D(2),p,\kappa\right)$ is SHD.

\end{cor}

In a similar vein to the question of $X$ versus $\beta X$, a natural question that arose in \cite{SHD} was whether it was possible to obtain something similar regarding $X$ versus $\mathscr{F}[X]$, the Pixley-Roy hyperspace of $X$. A. Bella and S. Spadaro proved in \cite{bella} the following theorem:

\begin{thm}\label{thm_tomita} If $X$ is SHD, then $\mathscr{F}[X]$ is SHD.
\end{thm}
Here, Theorem \ref{thm_tomita} is cited slightly different from \cite{bella}, by the fact that in the proof of Bella and Spadaro they don't use the property of being $T_1$, so, it can be omitted. 



On the other hand, similarly to what happened with $\beta X$, the presence of the SHD property in $\mathscr{F}[X]$ does not imply that $X$ has it. To establish this fact we first need to prove an auxiliary lemma. For what follows we will use the symbol $\CO(X)$ to denote the set formed by all clopen subsets of a topological space $X$.

\begin{lem}\label{lemma_jimenez_rios_tomita} Let $\kappa$ be an uncountable cardinal, $p$ be an element of $D(2)$ and $X$ be a dense subspace of $\Sigma(D(2),p,\kappa)$. If $\{U_n : n\in \omega\}$ is a subset of $\tau_{X}^{+}$ and $F\in \mathscr{F}[X]$, then there are $V\in \CO\left(D(2)^{\kappa}\right)$ and $\{x_n : n\in \omega\} \subseteq X$ such that $F\subseteq V$ and $x_n \in U_n \setminus V$ for every $n\in \omega$.

\end{lem}

\begin{pf} Let $m\in \omega\setminus 1$ and suppose that $\{y_i : 0\leq i \leq m\}$ is a faithful enumeration of $F$.

\medskip

\noindent {\bf Claim.} There exist $\{\alpha_i : 0\leq i \leq m\} \subseteq \kappa$ and $\{V(n,i) : n\in \omega \wedge 0 \leq i \leq m\} \subseteq \mathscr{B}\left(D(2)^{\kappa}\right)$ such that the following conditions hold for every $0\leq i \leq m$:

\begin{enumerate}
\item $\alpha_i \in \kappa \setminus \bigcup_{n\in \omega} \supp\left(V(n,i)\right)$;
\item $y_i \in \pi_{\alpha_i}^{-1}\{p\}$;
\item $V(n,0) \cap X \subseteq U_n$ for every $n\in \omega$ and
\item if $i>0$, then $V(n,i) \subseteq V(n,i-1) \setminus \pi_{\alpha_{i-1}}^{-1}\{p\}$.
\end{enumerate}

\medskip

We will do the construction by finite recursion. First, for every $n\in \omega$ there is $V(n,0) \in \mathscr{B}\left(D(2)^{\kappa}\right)$ such that $V(n,0) \cap X \subseteq U_n$. For this step we only need to notice that, since $y_0^{-1}\{p\}$ has $\kappa$ elements, there exists $\alpha_0 \in \kappa \setminus \bigcup_{n\in \omega} \supp\left(V(n, 0)\right)$ with $y_0(\alpha_0)=p$. 

Now suppose that for some $0\leq \ell < m$ we have produced $\{\alpha_i : 0\leq i \leq \ell\}$ and $\{V(n,i) : n\in \omega \wedge 0 \leq i \leq \ell\}$ with the desired properties. Since $\{V(n,\ell) \setminus \pi_{\alpha_\ell}^{-1}\{p\} : n\in \omega\}$ is formed by non-empty open subsets of $D(2)^{\kappa}$, for every $n\in \omega$ there is $V(n,\ell+1)\in \mathscr{B}\left(D(2)^{\kappa}\right)$ with $V(n,\ell+1)\subseteq V(n,\ell) \setminus \pi_{\alpha_\ell}^{-1}\{p\}$. Finally observe that, since $y_{\ell+1}^{-1}\{p\}$ has $\kappa$ elements, there is $\alpha_{\ell+1} \in \kappa \setminus \bigcup_{n\in \omega} \supp\left(V(n,\ell+1)\right)$ such that $y_{\ell+1}(\alpha_{\ell+1})=p$. This completes our construction.

To finish our proof set $V:= \bigcup_{i=0}^{m} \pi_{\alpha_i}^{-1}\{p\}$ and for every $n\in \omega$ let $x_n$ be an element of $X\cap \left(V(n,m) \setminus \pi_{\alpha_{m}}^{-1}\{p\}\right)$. Clearly, $V$ and $\{x_n : n\in \omega\}$ satisfy everything required.
\end{pf}

It turns out that for the class of spaces mentioned in Lemma~\ref{lemma_jimenez_rios_tomita} the Pixley-Roy hyperspace is always SHD.

\begin{thm}\label{thm_jimenez_rojas_tomita} If $\kappa$ is an uncountable cardinal, $p$ is an element of $D(2)$ and $X$ is a dense subspace of $\Sigma(D(2),p,\kappa)$, then $\mathscr{F}[X]$ is SHD.

\end{thm}

\begin{pf} Let $Y:= D(2)^{\kappa}$ and $\{[F_n,U_n] : n\in \omega\}$ be a sequence in $\mathscr{F}[X]$. The proof of the following statement will be by recursion.

\medskip

\noindent {\bf Claim 1.} There are $\{V_m : m\in \omega\} \subseteq \CO(Y)$ and $\{x(n,m) : (n,m)\in \omega\times \omega \} \subseteq X$ such that the following conditions hold for every $m\in \omega$:

\begin{enumerate}
\item $x(n,m)\in U_n\setminus V_m$ whenever $n\in \omega$ and
\item $F_m \cup \{x(m,i) : 0\leq i \leq m\}\subseteq V_{m+1}$.
\end{enumerate}

\medskip

For the base case, take any $q\in X$ and apply Lemma~\ref{lemma_jimenez_rios_tomita} to the sequence $\{U_n : n\in \omega\}$ and to the set $\{q\}$ to obtain $V_0\in \CO(Y)$ and $\{x(n,0) : n\in \omega\} \subseteq X$ such that $q \in V_0$ and $x(n,0)\in U_n\setminus V_0$ whenever $n\in \omega$. Now, suppose that for some $m\in \omega$ we have constructed $\{V_i : 0\leq i \leq m\}$ and $\{x(n,i) : n\in \omega \ \text{ and} \ 0 \leq i \leq m\}$ with the desired properties. By Lemma~\ref{lemma_jimenez_rios_tomita} applied to the sequence $\{U_n \setminus V_m : n\in \omega\}$ and to the set $F_m \cup \{x(m,i) : 0\leq i \leq m\}$ there exist $V_{m+1}\in \CO(Y)$ and $\{x(n,m+1) : n\in \omega\} \subseteq X$ such that $ F_m \cup \{x(m,i) : 0\leq i \leq m\} \subseteq V_{m+1}$ and $x(n,m+1)\in U_n\setminus V_{m+1 }$ whenever $n\in \omega$.

For each $m\in \omega$ let $G_m := F_m \cup \{x(m,i) : 0\leq i \leq m\}$. Clearly, $G_m \in [F_m,U_m]$ for all $m\in \omega$. Our next statement essentially concludes the proof.

\medskip

\noindent {\bf Claim 2.} The family $\left\{\{G_n\} : n\in \omega\right\}$ is locally finite, i.e., for every $G\in \mathscr{F}[X]$ there exists $V\in \tau_X$ with $G\subseteq V$ and such that $\left\{n\in \omega : G_n \in [G,V]\right\}$ is finite.

\medskip

 If $G \in \mathscr{F}[X]$ and for every $n\in\omega$ holds that $G \not \subseteq G_n$ , then $V:= X$ satisfies that $\left\{ n\in \omega : G_n \in [G,V]\right\} = \emptyset$. When there is $m\in \omega$ with $G \subseteq G_m$ what happens is that, for every $n> m+1$, holds that $x(n,m+1)\in G_n \cap \left( U_n\setminus V_{m+1}\right)$ and, therefore, $V:= V_{m+1} \cap X$ satisfies $G_n \not \in [G,V]$; consequently, $\left\{n\in \omega : G_n \in [G,V]\right\} \subseteq \{i\in \omega : 0\leq i \leq m+1\}$.

Finally, since the sequence $(G_n)$ satisfies that $\left\{\{G_n\} : n\in \omega\right\}$ is locally finite, it follows that $(G_n)$ has no convergent subsequences. Hence, $\mathscr{F}[X]$ is SHD.
\end{pf}

By Remark~\ref{ex_sigma_producto}, it follows that $\Sigma(D(2),p,\kappa)$ is not SHD whenever $\kappa\geq \mathfrak{s}$ and $p\in D(2)$. Theorem~\ref{thm_jimenez_rojas_tomita} yields the following corollary.

\begin{cor}\label{cor_jimenez_rojas_tomita} If $\kappa\geq \mathfrak{s}$ and $p\in D(2)$, then $\Sigma\left(D(2),p,\kappa\right)$ is not SHD and $\mathscr{F}\left[\Sigma\left(D(2),p,\kappa\right)\right]$ is SHD.

\end{cor}

\section{More divergence properties \textsc{I}}\label{MDPI}

We say that $X$ is {\it weakly selectively highly divergent} (in short, {\it WSHD}) if for every $\{U_n : n\in\omega\}\subseteq \tau_X^+$ there exists $A\in [\omega]^{\omega}$ such that, for every $n\in A$, we can choose $x_n\in U_n$ with the property that the sequence $(x_n)_{n\in A}$ has no convergent subsequences.

If $\{U_n : n\in\omega\}\subseteq \tau_X^+$ and $x\in X$, we shall say that $(U_n)$ {\it converges to $x$} if for every open set $U$ with $x\in U$ there is $m\in \omega$ such that $U_n \subseteq U$ whenever $n\geq m$. We will say that $X$ is {\it openly highly divergent} (in short, {\it OHD}) if there are no convergent sequences formed by elements of $\tau_X^+$.

A {\it local $\pi$-basis} for $x$ in $X$ is a family $\mathscr{B}\subseteq \tau_X^+$ with the property that for any $U\in \tau_X$ with $x\in U$ there exists $B\in\mathscr{B}$ such that $B\subseteq U$. When $\mathscr{B}$ is a local $\pi$-basis for $x$ in $X$ and $x \in \bigcap\mathscr{B}$ we say that $\mathscr{B}$ is a { \it local base} for $x$ in $X$.

The {\it local character} and {\it local $\pi$-character} of $x$ in $X$ are, respectively, the cardinal numbers \begin{align*} \chi(x,X) &:= \min\left\{\abs{\mathscr{B}} : \mathscr{B} \ \text{is a local basis for $x$ in $X$}\right\} \ \text{and} \\
\pi\chi(x,X) &:= \min\left\{\abs{\mathscr{B}} : \mathscr{B} \ \text{is a local $\pi$-basis for $x$ in $X$}\right\}.
\end{align*}

For more information about these cardinal functions and other ones, the reader can see \cite{hodel1984}. Finally, we will say that $X$ is a {\it UC} space (by uncountable character)(resp., {\it $\pi$-UC} (by uncountable $\pi$-character)) if for every $x\in X$ happens that $\chi(x,X)>\omega$ (resp., $\pi\chi(x,X)>\omega$).

Let us start with some basic relationships between these new concepts. It is clear that SHD implies WSHD, $\pi$-UC implies OHD and OHD implies UC. The fact that WSHD implies OHD requires a brief argument. Suppose that $X$ is not OHD and let $\{U_n : n\in \omega\}$ be a subset of $\tau_X^+$ that is convergent to $x\in X$. Then, if $A\in [\omega]^{\omega}$ and for every $n\in A$ we pick $x_n \in U_n$, it turns out that $(x_n)_{n\in A}$ converges to $x$. Thus, $X$ is not WSHD.

Unlike the SHD property (see Corollary~\ref{cor_thm_rios_tomita}), the OHD property {\bf is} inherited to dense subsets in the realm of regular spaces.

\begin{pro}\label{prop_OHD_denso} If $X$ is a regular OHD space and $D$ is a dense subset of $X$, then $D$ is OHD.

\end{pro}

\begin{pf} Let $\{U_n : n\in \omega\}$ be a subset of $\tau_D^+$ such that converges to $x\in D$. For each $n\in \omega$ let $V_n\in \tau_X^+$ satisfy $U_n = V_n \cap D$. Our goal is to prove that $(V_n)$ converges to $x$ in $X$. Let $V\in \tau_X$ be such that $x\in V$ and fix $W\in \tau_X$ with $x\in W \subseteq \cl_X W \subseteq V$. Observe that if $n\in \omega$ satisfies $U_n \subseteq W \cap D$, then $V_n \cap D \subseteq W\cap D$ and, therefore, $\cl_X \left(V_n \cap D \right) \subseteq \cl_X \left(W \cap D\right)$; that is, $\cl_X V_n \subseteq \cl_X W$, which in turn implies that $V_n \subseteq V$. Consequently, $\abs{\left\{n\in \omega : V_n \not \subseteq V\right\}} \leq \abs{\left\{n\in \omega : U_n \not \subseteq W \cap D\right\}}<\omega$.
\end{pf}

A standard fact is that if $p\in X$ and $\kappa$ is an infinite cardinal, then $\sigma\left(X,p,\kappa\right)$ is contained in $\Sigma\left(X,p,\kappa\right)$ and is dense in $X^{\kappa}$. This together with Proposition~\ref{prop_OHD_denso} implies the following result:

\begin{pro}\label{prop_OHD_potencia_1} Let $X$ be a regular space, $p\in X$ and $\kappa$ be an infinite cardinal. If $X^{\kappa}$ is OHD, then $\Sigma\left(X,p,\kappa\right)$ and $\sigma\left(X,p,\kappa\right)$ are OHD.

\end{pro}

On the other hand, if $X$ admits a sequence of non-empty open sets that is convergent to $p\in X$, and $p^{\omega} \in X^{\omega}$ denotes the sequence $p^{\omega}(n) := p$, then $ \sigma\left(X,p,\omega\right)$ admits a sequence of non-empty open sets that is convergent to $p^{\omega}$.

\begin{pro}\label{prop_OHD_potencia_2} If there exists a sequence $\{U_n : n\in \omega\}\subseteq \tau_X^{+}$ converging to $p\in X$, then there exists a sequence $\{V_n : n\in \omega\} \subseteq \tau_{\sigma\left(X,p,\omega\right)}^{+}$ converging to $p^{\omega}$.

\end{pro}

\begin{pf} For each $n\in \omega$ let $$V_n := \sigma\left(X,p,\omega\right)\cap \bigcap _{k=0}^{n} \pi _{k}^{- 1}[U_n].$$ If $W$ is an open subset of $\sigma\left(X,p,\omega\right)$ with $p^{\omega} \in W$, there exist $F\in [\omega]^{<\omega}\setminus\{\emptyset\}$ and $\left\{W_k : k\in F\right\} \subseteq \tau_X$ such that $$p^{ \omega}\in \sigma\left(X,p,\omega\right)\cap \bigcap _{k\in F} \pi _{k}^{-1}[W_k] \subseteq W.$$ Then, the condition $p\in \bigcap_{k\in F} W_k$ implies the existence of $m \in \omega$ such that $U_n \subseteq \bigcap_{k\in F} W_k$ for any $n\geq m $. Finally, a routine argument shows that $V_n \subseteq W$ for all $n\geq m + \max F$.
\end{pf}

Furthermore, if $\kappa$ is uncountable it is satisfied that the $\kappa$-th power of a $T_1$ space is OHD.

\begin{pro}\label{prop_OHD_potencia_3} If $X$ is a $T_1$ space and $\kappa$ is an uncountable cardinal, then $X^\kappa$ is OHD.

\end{pro}

\begin{pf} Let $x\in X^{\kappa}$ and $\{U_n : n\in \omega\}$ be a subset of $\mathscr{B}\left(X^{\kappa}\right)$. Since the set $\bigcup_{n\in \omega} \supp(U_n)$ is countable, we can take $\beta \in \kappa \setminus \bigcup_{n\in \omega} \supp(U_n)$. Let $U\in \tau_X$ satisfy $x(\beta) \in U$ and $X\setminus U \neq \emptyset$. Fix $n\in \omega$ and for every $\alpha\in \supp(U_n)$ let $y_\alpha \in \pi_\alpha[U_n]$ . Therefore, if $z\in X\setminus U$, then $$y := \left\{(\alpha,y_\alpha) : \alpha\in \supp(U_n)\right\} \cup \left(\left(\kappa\setminus \supp(U_n) \right)\times \{z\}\right)$$ is an element of $U_n \setminus \pi_{\beta}^{-1}[U]$, in particular, $U_n \not \subseteq \pi_{\beta}^{-1}[U]$. Thus, $(U_n)$ does not converge to $x$.
\end{pf}

Propositions~\ref{prop_OHD_potencia_1}, \ref{prop_OHD_potencia_2} and \ref{prop_OHD_potencia_3} are the proof of the following result.

\begin{thm}\label{thm_potencia} If $X$ is a $T_3$ space, $X$ admits a non-empty open sequence converging to $p\in X$ and $\kappa$ is an infinite cardinal, then the following statements are equivalent.

\begin{enumerate}
\item $X^{\kappa}$ is OHD.
\item $\Sigma\left(X,p,\kappa\right)$ is OHD.
\item $\sigma\left(X,p,\kappa\right)$ is OHD.
\item $\kappa>\omega$.
\end{enumerate}

\end{thm}

In particular, since the discrete space $D(2)$ satisfies the conditions of Theorem~\ref{thm_potencia}, we obtain the next corollary.

\begin{cor}\label{cor_OHD_Cantor} The following statements are equivalent for an infinite cardinal $\kappa$ and $p\in D(2)$.

\begin{enumerate}
\item $D(2)^{\kappa}$ is OHD.
\item $\Sigma\left(D(2),p,\kappa\right)$ is OHD.
\item $\sigma\left(D(2),p,\kappa\right)$ is OHD.
\item $\kappa>\omega$.
\end{enumerate}

\end{cor}

Finally, since $\Sigma\left(D(2),p,\kappa\right)$ is not WSHD because it is sequentially compact, Corollary~\ref{cor_OHD_Cantor} yields the following result.

\begin{cor}\label{corolarioacitar} For every cardinal $\kappa>\omega$ and $p\in D(2)$ it is satisfied that $\Sigma\left(D(2),p,\kappa\right)$ is OHD and not WSHD.

\end{cor}

\section{More divergence properties \textsc{II}}\label{MDPII}

The purpose of this section is to study the behavior of the OHD property between $X$ and $\mathscr{F}[X]$. To begin with, the following lemma will be useful to characterize the OHD property in $\mathscr{F}[X]$ in terms of $X$ when the latter is $T_2$.

\begin{lem}\label{lema_OHD_F[X]} If $X$ is a Hausdorff space and $F\in \mathscr{F}[X]$, then the following statements are equivalent.

\begin{enumerate}
\item $X$ has countable character at every point of $F$.
\item $\mathscr{F}[X]$ has countable character in $F$.
\item $\mathscr{F}[X]$ has countable $\pi$-character in $F$.
\item There is a sequence $\{V_n : n\in \omega\} \subseteq \tau_{\mathscr{F}[X]}^{+}$ converging to $F$.
\item There are $\{F_n : n\in \omega\} \subseteq \mathscr{F}[X]$ and $\{U_n : n\in \omega\} \subseteq \tau_X^+$ such that $\left\{[F_n, U_n] : n\in \omega\right\}$ converges to $F$.

\end{enumerate}

\end{lem}

\begin{pf} The implication $(1)\to (2)$ and the equivalence $(2) \leftrightarrow (3)$ can be found in \cite[Proposition~1, p.~336]{sakai1983}. Furthermore, it is evident that $(2) \to (4)$ and $(4) \to (5)$. To check that $(5) \to (1)$, suppose that $\{F_n : n\in \omega\} \subseteq \mathscr{F}[X]$ and $\{U_n : n\in \omega\} \subseteq \tau_X^+$ are such that $\left\{[F_n, U_n] : n\in \omega\right\}$ converges to $F$. Let $m\in \omega$ be such that if $n\geq m$, then $[F_n,U_n] \subseteq [F,X]$. Use the Hausdorff property to produce a family $\{V_x : x\in F\} \subseteq \tau_X^+$ such that $x\in V_x$ and $V_x \cap V_y = \emptyset$ whenever $x,y\in F$ are different.

Let $x$ be an element of $F$. Our goal is to show that $\{U_n \cap V_x : n\geq m\}$ is a local basis for $x$ in $X$. First, for any $n\geq m$ the inclusion $[F_n,U_n] \subseteq [F,X]$ implies that $F\subseteq U_n$; in particular, $x\in U_n$. On the other hand, if $U$ belongs to $\tau_X$ and $x\in U$, then for $V:= U \cup \left(\bigcup_{y\in F\setminus\{x\}} V_y\right)$ there is $k\in \omega$ such that $[F_n,U_n] \subseteq [F,V]$ when $n\geq k$. Finally, since $U_{m+k} \subseteq V$, it is satisfied that $U_{m+k}\cap V_x \subseteq V \cap V_x = U \cap V_x \subseteq U$.
\end{pf}

Our next theorem can be deduced from Lemma~\ref{lema_OHD_F[X]}.

\begin{thm}\label{thm_OHD_Pixley} The following statements are equivalent for a Hausdorff space $X$.

\begin{enumerate}
\item $\mathscr{F}[X]$ is OHD.
\item $\mathscr{F}[X]$ is UC.
\item $\mathscr{F}[X]$ is $\pi$-UC.
\item $X$ is UC.
\end{enumerate}

\end{thm}

By Theorem~\ref{thm_tomita} follows that any SHD space $X$ satisfies that $\mathscr{F}[X]$ is SHD. In fact, a simple modification of that argument can be used to prove that if $X$ is WSHD, then $\mathscr{F}[X]$ is WSHD. On the other hand, Theorem~\ref{thm_OHD_Pixley} guarantees that $\mathscr{F}[X]$ is OHD whenever $X$ is a Hausdorff UC space. What follows is to show that $\mathscr{F}[X]$ is WSHD if $X$ is regular and OHD. In order to do this we first prove the following lemma.

\begin{lem}\label{lema_tomita_rojas} Let $X$ be a regular OHD space. If $A$ is an infinite subset of $\omega$, $\{U_n : n\in A\}$ is a subset of $\tau_X^+$ and $F$ is an element of $\mathscr{F}[X]$, then there is a closed neighbourhood $V$ of $F$ (i.e., $V$ is a closed subset of $X$ and $F\subseteq \int V$) in such a way that the set $\left\{n\in A : U_n\setminus V \neq \emptyset\right\}$ is infinite.

\end{lem}

\begin{pf} Let $\{x_0,x_1,\ldots,x_m\}$ be an enumeration with no repetitions of $F$. Since $X$ is a regular OHD space, a finite recursion argument can be employed to construct a sequence $\{V_i : 0\leq i \leq m\} \subseteq P(X)$ and a decreasing sequence $\{A_i : 0\leq i \leq m\} \subseteq [A]^{\omega}$ such that: \begin{enumerate}
\item $V_i$ is a closed neighbourhood of $x_i$ for every $0\leq i \leq m$;
\item $A_0 \subseteq \left\{n\in A : U_n\setminus V_0 \neq \emptyset\right\}$; and 
\item $$A_i \subseteq \left\{n\in A_{i-1} : U_n\setminus \bigcup_{j=0}^{i} V_j \neq \emptyset\right\}$$ whenever $0<i\leq m$.
\end{enumerate} Hence, $V:= \bigcup_{i=0}^{m} V_i$ is a closed neighbourhood of $F$ and $A_m$ is an infinite subset of $\left\{n\in A : U_n\setminus V \neq \emptyset\right\}$.
\end{pf}

With the previous result available, we are all set to prove one of the central theorems of this section.

\begin{thm}\label{thm_rojas_tomita} If $X$ is a regular OHD space, then $\mathscr{F}[X]$ is WSHD.

\end{thm}

\begin{pf} Let $\{F_n : n\in \omega\} \subseteq \mathscr{F}[X]$ and $\{U_n : n\in \omega\} \subseteq \tau_X$ be such that $\left\{[F_n,U_n] : n\in \omega\right\}$ is a sequence of non-empty open subsets of $\mathscr{F}[X]$. The following statement is the cornerstone of our proof.

\medskip

\noindent {\bf Claim.} There are a decreasing sequence $\{A_k : k\in \omega\} \subseteq [\omega]^{\omega}$, $\{G_k : k\in \omega\} \subseteq \mathscr{F}[X]$, $\{V_k : k\in \omega\} \subseteq P(X)$ and $\{x(n,k) : (n,k)\in \omega\times \omega\}$ in such a way that the following properties hold for every $k\in \omega$:

\begin{enumerate}
\item $V_0$ is a closed neighbourhood of $F_0$;
\item for each $n\in A_k$, $x(n,k)\in (U_n\setminus F_n)\setminus \bigcup_{i=0}^{k} V_i$;
\item $G_k := F_k \cup \left\{x(k,m) : m\leq k  \wedge k\in A_m\right\}$; and
\item $V_k$ is a closed neighbourhood of $\bigcup_{i=0}^{k-1} G_i$ whenever $k>0$.
\end{enumerate}

\medskip

We proceed by recursion. Since our hypothesis implies that $\{U_n\setminus F_n : n\in \omega\}$ is a sequence of non-empty open subsets of $X$, Lemma~\ref{lema_tomita_rojas} produces $A_0\in [\omega]^{\omega}$ and $V_0$, a closed neighbourhood of $F_0$, such that $(U_n\setminus F_n) \setminus V_0 \neq \emptyset$ for every $n\in A_0$. For each $n\in A_0$ let $x(n,0)$ be an element of $(U_n\setminus F_n) \setminus V_0$. Set $$G_0 := \begin{cases} F_0\cup \{x(0,0)\}, & \text{if} \ 0\in A_0, \\
F_0, & \text{if} \ 0\not\in A_0,
  \end{cases}$$ i.e., $G_0 = F_0\cup \left\{x(0,m) : m\leq 0  \wedge 0\in A_0\right\}$.
  
Now suppose that for $k\in \omega$ we have constructed $\{A_i : 0\leq i \leq k\}$, $\{G_i : 0\leq i \leq k\}$, $\{V_i : 0\leq i \leq k\}$ and $\{x(n,i) : 0\leq i \leq k \wedge n\in A_i\}$ with the desired properties. From Lemma~\ref{lema_tomita_rojas} we can obtain $A_{k+1}\in [A_k]^{\omega}$ and $V_{k+1}$, a closed neighbourhood of $\bigcup_{i=0}^{k} G_i$, such that $(U_n\setminus F_n)\setminus \bigcup_{i=0}^{k+1} V_i\neq\emptyset$ for every $n\in A_{k+1}$.  For each $n\in A_{k+1}$ let $x(n,k+1)$ be an element of $(U_n\setminus F_n)\setminus \bigcup_{i=0}^{k+1} V_i$. Finally, set $$G_{k+1} := F_{k+1}\cup \left\{x(k+1,m) : m\leq k+1  \wedge k+1\in A_m\right\}.$$ This completes our construction.

It is clear that $G_n\in [F_n,U_n]$ for every $n\in \omega$. Now, recursively generate a strictly increasing sequence $(n_k)$ in such a way that $n_k \in A_k$ for every $k\in \omega$. Our final goal is to prove that, if $A: = \{n_k : k\in \omega\}$, then $(G_n)_{n\in A}$ has no convergent subsequences.

Before we begin we need to make a key remark: if $j\in \omega$ and $n\in A$ satisfy $n_j \leq n$, then $n\in A_j$ and $x(n,j) \in G_n$. Indeed, if $k\in \omega$ is such that $n = n_k$, then $j\leq k$ and thus, the relations $n = n_k \in A_k \subseteq A_j$ imply $n\in A_j$; consequently, $x(n,j)\in G_n$ since $j\leq n_j\leq n$.

To finish our proof suppose, in search of a contradiction, that $B\in [A]^{\omega}$ and $H\in \mathscr{F}[X]$ are such that $(G_n)_{n\in B}$ converges to $H$ in $\mathscr{F}[X]$. Let $s \in \omega$ satisfy $\{G_n : n\in B \cap [s,\omega)\} \subseteq [H,X]$ and set $m := \min (B\setminus s)$. Notice that since $H\subseteq G_m \subseteq \int V_{m+1}$, then $[H, \int V_{m+1}]$ is an open neighbourhood of $H$ and therefore, $\{G_n : n\in B \cap [t,\omega)\} \subseteq [H, \int V_{m+1}]$ for some $t\in \omega$.

Let $k\in \omega$ satisfy $m= n_k$ and fix $j\in \omega$ with $n_k<j$. For every $n\in B$ such that $n_j \leq n$, our key remark ensures that $x(n,j) \in G_n$. Furthermore, since $x(n,j)\not\in \bigcup_{i=0}^{j} V_i$ and $m = n_k < j$, it follows that $x(n,j) \not \in V_{m+1}$; in particular, $x(n,j) \in G_n \setminus \int V_{m+1}$. Consequently, $B\cap [n_j,\omega)$ is an infinite subset of $\left\{n\in B : G_n \not \in [H, \int V_{m+1}]\right\}$, which contradicts the last line of the previous paragraph.
\end{pf}

Recall that $X$ is called a $P$-space if the countable intersection of open sets is also open. We will show next that in the realm of regular $P$-spaces and Hausdorff topological groups a lot of the properties we have presented so far are, in fact, equivalent.

\begin{pro}\label{prop_OHD_P-space} If $X$ is a regular $P$-space and $x\in X$, then the following statements are equivalent.

\begin{enumerate}
\item $X$ has countable character in $x$.
\item There exists $\{U_n : n\in \omega\}\subseteq \tau_X^+$ converging to $x$.
\item $X$ has countable $\pi$-character in $x$.
\end{enumerate}

\end{pro}

\begin{pf} We only need to see that $(3) \to (1)$. Suppose that $\{V_n : n\in \omega\}\subseteq \tau_X^+$ is a local $\pi$-base for $x$ in $X$. Since every regular $P$-space is zero-dimensional (see \cite[1W, p.~69]{Porter}), for every $n\in \omega$ there exists $C_n$, a non-empty clopen subset of $X$, with $C_n \subseteq V_n$. Evidently, $\{C_n : n\in \omega\}$ is a local $\pi$-basis for $x$ in $X$.

\medskip

\noindent {\bf Claim.} If $I := \left\{n\in \omega : x\in C_n\right\}$, then $\{C_n : n\in I\}$ is a local basis for $X $ in $x$.

\medskip

Indeed, if $U\in \tau_X$ satisfies $x\in U$ and $V:= U\cap \bigcap_{n\in \omega\setminus I} (X\setminus C_n)$, then $V$ is an open subset of $X$ with $x\in V$. Hence, there exists $n\in \omega$ with $C_n \subseteq V$ and thus, the inclusion $C_n \subseteq \bigcap_{m\in \omega\setminus I} (X\setminus C_m)$ implies that $n\in I$; consequently, $x\in C_n \subseteq U$.
\end{pf}

Thus, we have the following result.

\begin{thm}\label{thm_P-space} If $X$ is a regular $P$-space, then the following statements are equivalent.

\begin{enumerate}
\item $\mathscr{F}[X]$ is SHD.
\item $\mathscr{F}[X]$ is WSHD.
\item $\mathscr{F}[X]$ is OHD.
\item $\mathscr{F}[X]$ is UC.
\item $\mathscr{F}[X]$ is $\pi$-UC.
\item $X$ is SHD.
\item $X$ is WSHD.
\item $X$ is OHD.
\item $X$ is UC.
\item $X$ is $\pi$-UC.
\end{enumerate}

\end{thm}

\begin{pf} The implications $(1)\to (2) \to (3) \to (4)$ and $(6)\to (7) \to (8) \to (9)$ always hold. Since $X$ is a Hausdorff space, Theorem~\ref{thm_OHD_Pixley} ensures that (3) -- (5) and (9) are equivalent. On the other hand, the fact that (8) -- (10) are equivalent follows from Proposition~\ref{prop_OHD_P-space}. Moreover, $(6) \to (1)$ by virtue of Theorem~\ref{thm_tomita}. Notice that a straightforward argument shows that if $X$ is a $T_1$ and UC $P$-space then it satisfies that every open set is infinite and $X$ contains no non-trivial convergent sequences and thus, implication $(9) \to (6)$ is due to \cite[Proposition~1, p.~4]{SHD}.
\end{pf}

\begin{pro}\label{prop_OHD_grupo_top} If $G$ is a topological group, $x\in G$ and $e\in G$ is the neutral element of $G$, then the following statements are equivalent.

\begin{enumerate}
\item $G$ has countable character in $e$.
\item $G$ has countable character in $x$.
\item There exists $\{U_n : n\in \omega\}\subseteq \tau_G^+$ converging to $x$.
\item $G$ has countable $\pi$-character in $x$.
\item $G$ has countable $\pi$-character in $e$.
\end{enumerate}

\end{pro}

\begin{pf} Clearly,  $(2) \to (3) \to (4)$. On the other hand, the homogeneity of $G$ certifies that $(1)\to (2)$ and $(4)\to (5)$. To see that $(5)\to (1)$ suppose that $\{V_n : n\in \omega\}\subseteq \tau_G^+$ is a $\pi$-local basis for $G$ in $e$. Our goal is to verify that $\left\{V_n \cdot \left(V_n\right)^{-1} : n\in \omega\right\}$ is a local basis for $G$ in $e$. Indeed, if $U\in \tau_G$ satisfies $e\in U$, then there exists $V\in \tau_G$ with $e\in V$ and $V\cdot V^{-1} \subseteq U$. Then, if $n\in \omega$ satisfies $V_n \subseteq V$, it follows that $e\in V_n \cdot \left(V_n\right)^{-1} \subseteq U$.
\end{pf}

Theorem~\ref{thm_OHD_Pixley} and Proposition~\ref{prop_OHD_grupo_top} allow us to obtain a result similar to Theorem~\ref{thm_P-space} in the realm of Hausdorff topological groups.

\begin{thm}\label{thm_top_gp} If $G$ is a Hausdorff topological group, then the following statements are equivalent.

\begin{enumerate}
\item $\mathscr{F}[G]$ is WSHD.
\item $\mathscr{F}[G]$ is OHD.
\item $\mathscr{F}[G]$ is UC.
\item $\mathscr{F}[G]$ is $\pi$-UC.
\item $G$ is OHD.
\item $G$ is UC.
\item $\chi(e,G)>\omega$.
\item $G$ is $\pi$-UC.
\item $\pi\chi(e,G)>\omega$.
\end{enumerate}

\end{thm}

\begin{pf} Since $G$ is a Hausdorff space, Theorem~\ref{thm_OHD_Pixley} guarantees that (2) -- (4) and (6) are equivalent. On the other hand, (5) -- (9) are equivalent thanks to Proposition~\ref{prop_OHD_grupo_top}. To finish the proof, recall that $(1)\to (2)$ is always true and notice that $(5)\to (1)$ by virtue of Theorem~\ref{thm_rojas_tomita}.
\end{pf}

Unfortunately, it is not possible to add the statement {\lq\lq}$G$ is WSHD{\rq\rq} in the list of Theorem~\ref{thm_top_gp} since $\Sigma\left(D(2),0 ,\mathfrak{c}\right)$ is a Hausdorff topological group that is not WSHD by Corollary~\ref{corolarioacitar} and its Pixley-Roy hyperspace is SHD and therefore WSHD by Corollary~\ref{cor_jimenez_rojas_tomita}.

Finally, if $X$ is a Tychonoff space and $C_p(X)$ stands for the ring of continuous function equipped with the topology of pointwise convergence, it is well-known that $C_p(X)$ is a topological group. In this particular case we can prove the following strengthening of Theorem~\ref{thm_top_gp}.

\begin{thm}\label{thm_jimenez_rios_tomita} The following statements are equivalent for a Tychonoff space $X$.

\begin{enumerate}
\item $\mathscr{F}\left[C_p(X)\right]$ is SHD.
\item $\mathscr{F}\left[C_p(X)\right]$ is WSHD.
\item $\mathscr{F}\left[C_p(X)\right]$ is OHD.
\item $\mathscr{F}\left[C_p(X)\right]$ is UC.
\item $\mathscr{F}\left[C_p(X)\right]$ is $\pi$-UC.
\item $C_p(X)$ is SHD.
\item $C_p(X)$ is WSHD.
\item $C_p(X)$ is OHD.
\item $C_p(X)$ is UC.
\item $\chi\left(\overline{0},C_p(X)\right)>\omega$.
\item $C_p(X)$ is $\pi$-UC.
\item $\pi\chi\left(\overline{0},C_p(X)\right)>\omega$.
\item $X$ is uncountable.
\end{enumerate}

\end{thm}

\begin{pf} Items (2) -- (5) and (8) -- (12)  are all equivalent to each other by virtue of Theorem~\ref{thm_top_gp}. Furthermore, $(1)\to (2)$ and $(6)\to (7) \to (8)$ always hold. On the other hand, $(6) \to (1)$ follows from Theorem~\ref{thm_tomita} and $(10) \to (13)$ since $\chi\left(C_p(X)\right)=|X|$ (see \cite[Theorem~\textsc{I}.1.1, p.~25]{ arh1992}). The argument for $(13) \to (6)$ goes as follows. 

Let $\{U_n: n\in\omega \}$ be a subset of $\tau_{C_p(X)}^{+}$ formed by canonical elements of the standard base. For every $n \in \omega$, let $f_n\in C_p(X)$, $\varepsilon_n \in (0,\infty)$ and $F_n \in [X]^{<\omega}$ be such that $U_n = [f_n; F_n; \varepsilon_n]$. Notice that, since $X$ is uncountable, there is a point $z\in X\setminus \bigcup_{n\in\omega} F_n$. For each $n\in \omega$, let $G_n := F_n \cup \{z\}$ and let $g_n : G_n  \to \mathbb{R}$ be defined as $$g_n(x) := \begin{cases} f_n(x), & \text{if} \ x\in F_n, \\
n, & \text{if} \ x=z.
  \end{cases}$$ Since every $G_n$ is $C^{*}$-embedded in $X$, there is a continuous function $h_n : X \to \mathbb{R}$ such that its restriction to $G_n$ coincides with $g_n$; in particular, $h_n \in U_n$ for every $n\in \omega$. However, for every $h\in C_p(X)$ the open set $[h;\{z\};1/2]$ intersects $\{h_n : n\in \omega\}$ in at most one point; consequently, $(h_n)$ has no convergent subsequences.
\end{pf}

We have mentioned several times that there are no first countable SHD spaces. A natural question is whether there are SHD spaces that are Fréchet-Urysohn. Since for a compact Hausdorff space $X$ it is satisfied that $C_p(X)$ is Fréchet-Urysohn if and only if $X$ is scattered (see \cite[Theorem~\textsc{III}.1.2, p.~91]{arh1992}), Theorem~\ref{thm_jimenez_rios_tomita} provides a list of examples of the sort; namely, any uncountable, scattered, compact and Hausdorff space $X$ satisfies that $C_p(X)$ is a Fréchet-Urysohn SHD space. In particular, $C_p\left([0,\omega_1]\right)$ has the desired properties.

\section{More divergence properties \textsc{III}}\label{MDPIII}

The space $X$ will be called {\it sequentially discrete} (in short, {\it SD}) if every convergent sequence in $X$ is eventually constant. It was proven in \cite[Proposition~1, p.~4]{SHD} that any SD space in which every non-empty open set is infinite is SHD. On the other hand, a Tychonoff space $X$ is an {\it $F'$-space} if any two disjoint cozero sets have disjoint closures. It is known that any $F'$-space is SD (see \cite[Fact~3.2, p.~574]{hertam2011}).

Our first task in the present section is to prove that a $T_1$ space is SD if and only if its Pixley-Roy hyperspace is SD. We begin with an auxiliary result.

\begin{lem}\label{lema_Elmer} If $X$ a $T_1$ space and $\mathscr{F}[X]$ is not SD, then there is an injective sequence $\{F_n : n\in \omega\} \subseteq \mathscr{F}[X]$ in such a way that, if $F:= \bigcap_{n\in \omega} F_n$, then $(F_n)$ converges to $F$.

\end{lem}

\begin{pf} Let $\{G_n : n\in \omega\}$ be a non-trivial sequence in $\mathscr{F}[X]$ in such a way that $(G_n)$ converges to a set $G\in \mathscr{F}[X]$. Firstly, since $(G_n)$ is not eventually the constant $G$, the set $A:= \{n\in \omega : G_n \neq G\}$ is infinite. Notice that since $(G_n)_{n\in A}$ is a sequence converging to $G$, $G\not\in \{G_n : n\in A\}$ and $\mathscr{F}[X]$ is a $T_1$ space, it must be the case that $\{G_n : n\in A\}$ is an infinite set. Let $\{F_n : n\in \omega\}$ be an enumeration with no repetitions of $\{G_n : n\in A\}$.

Now, since $[G,X]$ is an open set in $\mathscr{F}[X]$, there is $m\in \omega$ such that $\{F_n : n\geq m\} \subseteq [G,X]$, which in turn implies that $G \subseteq \bigcap_{n\geq m} F_n$. Finally, if $F:= \bigcap_{n\geq m} F_n$ and $U\in \tau_X$ satisfies $F \subseteq U$, then $[G,U]$ is an open set in $\mathscr{F}[X]$ and therefore, there is $k\geq m$ in such a way that $F_n \in [G,U]$ for every $n\geq k$; consequently, $F_n \in [F,U]$ for whenever $n\geq k$. Thus, $\{F_n : n\geq m\}$ satisfies the properties required in our lemma.

\end{pf}

Equipped with the above lemma we are ready to prove that the SD property appears in $X$ if and only if it does so in $\mathscr{F}[X]$.

\begin{pro}\label{auxilio} If $X$ is a $T_1$ space, then $X$ is SD if and only if $\mathscr{F}[X]$ is SD.

\end{pro}

\begin{pf} Since every convergent sequence in $X$ naturally induces a convergent sequence in $\mathscr{F}[X]$, it follows that $X$ is SD whenever $\mathscr{F}[X]$ is SD. For the remaining implication let us assume by contradiction that $X$ is SD and $\mathscr{F}[X]$ is not SD. Use Lemma~\ref{lema_Elmer} to fix an injective sequence $\{F_n : n\in \omega\} \subseteq \mathscr{F}[X]$ in such a way that, if $F:= \bigcap_{n\in \omega} F_n$, then $(F_n)$ converges to $F$.

Let $A:= \{n\in \omega : F_n \neq F\}$, pick $x_n \in F_n\setminus F$ for every $n\in A$ and let $\{z_1,\ldots,z_k\}$ be a faithful enumeration of $F$. Since $X$ is a $T_1$ SD space, a finite recursion argument can be used to construct a sequence $\{U_i : 1\leq i \leq k\} \subseteq \tau_X^ +$ and a decreasing sequence $\{A_i : 1\leq i \leq k\} \subseteq [A]^{\omega}$ such that, for every $1\leq i \leq k$, $U_i \cap F = \{z_i\}$ and $$\{x_n : n\in A_i\} \cap \bigcup_{j=1}^{i} U_i = \emptyset.$$

Finally, if $U := \bigcup_{i=1}^k U_i$, then $[F,U]$ is an open set in $\mathscr{F}[X]$ and thus, there is $m\in \omega$ such that $F_n \in [F,U]$ for every $n\in A_k \cap [m,\infty)$; consequently, $\{x_n : n\in A_k \cap [m,\infty)\} \subseteq U$, a contradiction. Hence, $\mathscr{F}[X]$ is SD.

\end{pf}

Gillman constructed in \cite{EHD} a $P$-space $X$ and an extremally disconnected space $Y$ such that $X\times Y$ is Tychonoff and not an $F'$-space. Notice that $X\times Y$ is SD since $X$ and $Y$ are SD (see \cite[6L, p.~505]{Porter}). Hence, we have proven the following statement.

\begin{pro} There exists a Tychonoff SD space which is not an $F'$-space.

\end{pro}

On the other hand, there are examples of Pixley-Roy hyperspaces which are SHD but not SD.

\begin{pro} There exists a Hausdorff space $X$ such that $\mathscr{F}[X]$ is SHD and not SD.

\end{pro}

\begin{pf} Let $X$ be a Hausdorff SHD space with a non-trivial convergent sequence $(x_n)$ (see \cite[Example~2, p.~5]{SHD}). It follows from Theorem~\ref{thm_tomita} that $\mathscr{F}[X]$ is SHD and by Proposition~\ref{auxilio} that $\mathscr{F}[X]$ is not SD.
\end{pf}

Regarding the $F'$ property between $X$ and $\mathscr{F}[X]$, it turns out that this concept does not reflect from $X$ onto the hyperspace $\mathscr{F}[X]$.

\begin{pro} There exists a compact Hausdorff $F'$-space $X$ such that $\mathscr{F}[X]$ is not an $F'$-space.

\end{pro}

\begin{pf} Let $Y$ and $Z$ be disjoint copies of $\beta\omega$ such that $\beta\omega= Y\cup Z$. Now, \cite[Theorem~2.6, p.~145]{lutzer} implies that \begin{align*} \mathscr{F}[\beta\omega] &\cong  \mathscr{F}[Y]\oplus \mathscr{F}[Z]\oplus \left(\mathscr{F}[Y]\times \mathscr{F}[Z] \right)\\
&\cong  \mathscr{F}[\beta\omega]\oplus \mathscr{F}[\beta\omega]\oplus \left(\mathscr{F}[\beta\omega]\times \mathscr{F}[\beta\omega] \right).
\end{align*}

Looking for a contradiction, suppose that $\mathscr{F}[\beta\omega]$ is an $F'$-space. Since $\mathscr{F}[\beta\omega]\times\mathscr{F}[\beta\omega]$ is $C^{*}$-embedded in $\mathscr{F}[\beta\omega]\oplus \mathscr{F}[\beta\omega]\oplus \left(\mathscr{F}[\beta\omega]\times \mathscr{F}[\beta\omega] \right)$, it follows that $\mathscr{F}[\beta\omega]\times\mathscr{F}[\beta\omega]$ is an $F'$-space (see \cite[Theorem~0.1, p.~276]{adowfspaces}) and therefore, $\mathscr{F}[\beta\omega]$ is a $P$-space (see \cite[Theorem, p.~51]{curtis}). Then, it follows from \cite[Proposition~5]{SHD} that $\beta\omega$ is a $P$-space, which in turn implies that $\beta\omega$ is finite (\cite[4AG, p.~353]{Porter}). In sum, $\beta\omega$ is an $F'$-space and $\mathscr{F}[\beta\omega]$ is not. 
\end{pf}

We mentioned at the beginning of Section~\ref{MDPI} that $\pi$-UC implies OHD and OHD implies UC. We will show below that these implications are not reversible. For our next results the symbol $\omega^{*}$ will denote the space $\beta\omega\setminus \omega$, i.e., the remainder of the Stone-\v{C}ech compactification of $\omega$.

\begin{pro} There exists a UC space which is not OHD.

\end{pro}

\begin{pf} For each $\alpha \in \omega_1$ let $X_\alpha$ stand for the space $\omega^{*} \times \{\alpha\}$. Define $X:= \bigoplus_{\alpha\in\omega_1} X_\alpha$, fix a point $p\not \in X$ and set $Y:= X\cup \{p\}$. Declare $X$ to be an open subset of $Y$ and let $$\left\{Y\setminus \bigcup_{\alpha \in F} X_\alpha : F \in [\omega_1]^{<\omega}\right\}$$ be a local basis for $Y$ in $p$.

Clearly, each member of $Y\setminus \{p\}$ has uncountable character. Now, if $\{F_n : n\in\omega\}$ is a subset of $[\omega_1]^{<\omega}$ and $\beta \in \omega_1 \setminus \bigcup_{n\in\omega} F_n$, then for every $n\in \omega$ it is satisfied that $$Y\setminus \bigcup_{\alpha \in F_n} X_\alpha \not\subseteq Y\setminus X_\beta.$$ Thus, $\chi(p,Y)>\omega$. Finally, since $\{X_n : n\in \omega\}$ is a sequence in $\tau_Y^+$ converging to $p$, it follows that $Y$ is not OHD.
\end{pf}

Regarding OHD versus $\pi$-UC, we can even produce an SD space in which every open set is infinite and that is not $\pi$-UC. In fact, any SD space in which every open set is infinite is automatically a SHD space. 

\begin{pro} There exists an SD and OHD space in which every open set is infinite such that is not $\pi$-UC.

\end{pro}

\begin{pf} For each $n\in \omega$ let $X_n$ be the space $\omega^{*}\times\{n\}$. For every $U\in \tau_{\beta\omega}$ define $$U^{*}:=(U\cap \omega^{*})\cup \left(\bigcup_{n\in U} X_n\right).$$ Set $X:= \bigoplus_{n\in \omega} X_n$ and $Y:= X\cup \omega^{*}$. Declare $X$ to be an open subset of $Y$ and for every $p\in \omega^{*}$ let $\{U^{*} : p\in U \wedge U\in \tau_{\beta\omega}\}$ be a local basis for $Y$ in $p$. Clearly, every open set in $Y$ is infinite.

Evidently, $\{X_n : n\in\omega\}$ is a local $\pi$-base for $Y$ in every point of $\omega^{*}$; therefore, $Y$ is not $\pi$-UC. In order to prove that $Y$ is SD, let $(y_n)$ be an injective sequence in $Y$. Since each $X_n$ is open in $Y$ and is homeomorphic to the SD space $\omega^{*}$, it cannot happen that $(y_n)$ converges to a point in $X$. Now, if $A := \{n\in \omega: y_n\in\omega^{*} \}$ is infinite, then $(y_n)_{n\in A}$ is a non-trivial sequence in the SD space $\omega^{*}$ and thus, $(y_n)$ cannot converge to any element of $\omega^{*}$.

In the case where $A$ is finite, for every $n\in \omega\setminus A$ let $m(n) \in \omega$ satisfy $y_n \in X_{m(n)}$ and fix $y\in \omega^{*}$. Since the sequence $(m(n))_{n\in\omega\setminus A}\subseteq\omega$ does not converge to $y$ in $\beta\omega$, there is $U\in \tau_{\beta\omega}$ with $y\in U$ such that $B:=\{n\in\omega\setminus A : m(n)\not\in U\}$ is infinite. Notice that for every $n\in B$ the condition $m(n) \not\in U$ implies that $U^{*}\cap X_{m(n)}=\emptyset$; in particular, $y_n \not \in U^{*}$. Thus, $\{y_n\}_{n\in B}$ is disjoint from $U^*$ and hence, $\{y_n \}_{n\in \omega}$ does not converge to $y$ in $Y$.

\end{pf}

To conclude this section we include a diagram that summarizes most of the results presented throughout the text. It is necessary to mention that the implications marked with (*) need non-empty open sets to be infinite in the corresponding space to be true. Finally, the indication {\lq\lq}$+T_i${\rq\rq} means that $X$ needs to be a $T_i$ space for the implication to be valid.

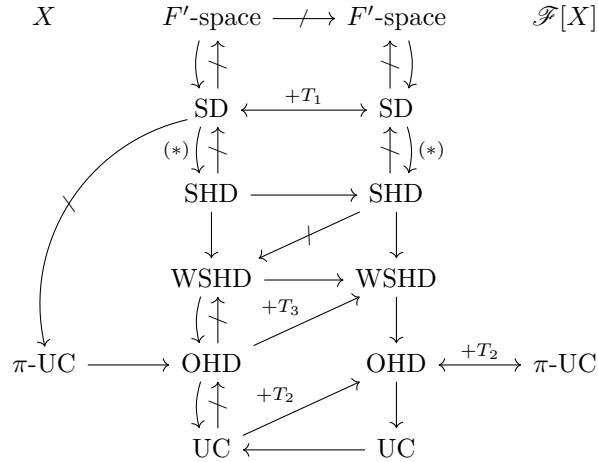
\begin{center}

\begin{tikzcd}

X & F'\text{-space} \arrow[d, shift right, bend right = 15]{}{} \arrow[r, "/"{anchor=center,sloped}]{}{} & F'\text{-space} \arrow[d, shift left, bend left = 15]{}{} & \mathscr{F}[X] \\

& \text{SD} \arrow[dddl, shift left, "/"{anchor=center,sloped}, bend right = 45]{}{} \arrow[d, swap, shift right, bend right = 15]{}{(*)} \arrow[r, leftrightarrow]{}{+T_1} \arrow[u, shift right, "/"{anchor=center,sloped}]{}{} & \text{SD} \arrow[d, shift left, bend left = 15]{}{(*)} \arrow[u, shift left, "/"{anchor=center,sloped}]{}{} & \\

& \text{SHD} \arrow[d]{}{} \arrow[r]{}{} \arrow[u, shift right, "/"{anchor=center,sloped}]{}{} & \text{SHD} \arrow[d]{}{} \arrow[u, shift left, "/"{anchor=center,sloped}]{}{} \arrow[ld, "/"{anchor=center,sloped}]{}{} & \\

& \text{WSHD} \arrow[d, shift right, bend right = 15]{}{} \arrow[r]{}{} & \text{WSHD} \arrow[d]{}{} & \\

\pi\text{-UC} \arrow[r]{}{} & \text{OHD} \arrow[d, shift right, bend right = 15]{}{} \arrow[ru]{}{+T_3} \arrow[u, shift right, "/"{anchor=center,sloped}]{}{} & \text{OHD} \arrow[d]{}{} \arrow[r, leftrightarrow]{}{+T_2} & \pi\text{-UC} \\

& \text{UC} \arrow[ru]{}{+T_2} \arrow[u, shift right, "/"{anchor=center,sloped}]{}{} & \text{UC} \arrow[l]{}{} & \\

\end{tikzcd}

\end{center}

\section{Questions}

In this article we have answered Question~1, Question~2 in an independent way from Bella and Spadaro's answer in \cite{bella}, Question~6 and Question~7 from \cite{SHD}. Furthermore, new questions have arisen from the topics presented in this text.

\begin{question}\label{Q1} Is there a Hausdorff zero-dimensional SHD space $X$ such that $\beta X$ is not zero-dimensional?

\end{question}

The natural candidate to solve Question~\ref{Q1} is the space constructed in \cite{kph} but we have not been able to determine if it does not have isolated points.

\begin{question}\label{Q2} Does it hold that if $X$ is a Tychonoff non-compact SHD space, then $\beta X$ is SHD?

\end{question}

Theorem~\ref{thm_rios_tomita} reinforces the hope of having a positive answer to Question~\ref{Q2}.

\begin{question}\label{Q3} Is it true that if $X$ is Tychonoff, non-compact, locally compact and Lindelöf, then $\beta X\setminus X$ is SHD?

\end{question}

Regarding Question~\ref{Q3}, the statement obtained by changing the word {\lq\lq}Lindelöf{\rq\rq} to {\lq\lq}$\sigma$-compact{\rq\rq} was proved in \cite{SHD}.

Finally, concerning Sections~\ref{MDPI}, \ref{MDPII} and \ref{MDPIII} of this article, the following questions remain unanswered.

\begin{question} Does it hold that every WSHD space is SHD?

\end{question}

\begin{question} Is $X$ OHD whenever $\mathscr{F}[X]$ is WSHD?

\end{question}

\begin{question} Are the following statements equivalent?

\begin{enumerate}
\item $\mathscr{F}[X]$ is SHD.
\item $\mathscr{F}[X]$ is WSHD.
\item $\mathscr{F}[X]$ is OHD.

\end{enumerate}

\end{question}

\begin{question} Is it true that if $\mathscr{F}[X]$ is an $F'$-space, then $X$ is an $F'$-space?

\end{question}

\end{document}